\documentclass{article}[12pt]
\usepackage{amsmath,amssymb}

\title{Oracle-supported  drawing of the Gr\"obner {\em escalier}}

\author{
{\small\bf Maria Emilia Alonso} \\
{\small Depto. Algebra. Fac.CC. Matem\'aticas}\\
{\small Universidad Complutense de Madrid}\\
{\small{\tt m\_alonso@mat.ucm.es}}
\and
{\small\bf Maria Grazia Marinari}\\
{\small DIMA}\\
{\small Universit\`a di Genova}\\
{\small{\tt marinari@dima.unige.it}}
\and
{\small\bf Teo Mora}\\
{\small DISI}\\
{\small Universit\`a di Genova}\\
{\small{\tt theomora@disi.unige.it}}
}

\def\And{\mbox{\ and }}

\def\Bbb#1{{\mathbb #1}}
\def\Cal#1{{\cal #1}}

\def\Can{\mathop{\rm Can}\nolimits}
\def\then{\;\Longrightarrow\;}

\def\Span{\mathop{\rm Span}\nolimits}
\def\lc{\mathop{\rm lc}\nolimits}
\def\supp{\mathop{\rm supp}\nolimits}

\def\lcm{\mathop{\rm lcm}\nolimits}

\begin{document}
\maketitle
\def\qed{\ifmmode\squareforqed\else{\unskip\nobreak\hfil
\penalty50\hskip1em\null\nobreak\hfil\squareforqed
\parfillskip=0pt\finalhyphendemerits=0\endgraf}\fi}
\def\squareforqed{\hbox{\rlap{$\sqcap$}$\sqcup$}}
\def\pf{\par\noindent {\it Proof }}
\newtheorem{Theorem}{Theorem}
\newtheorem{Corollary}[Theorem]{Corollary}
\newtheorem{Lemma}[Theorem]{Lemma}
\newtheorem{Proposition}[Theorem]{Proposition}
\newtheorem{Fact}[Theorem]{Fact}
\newtheorem{Definition}[Theorem]{Definition}
\newtheorem{Problem}[Theorem]{Problem}
\newtheorem{Example}[Theorem]{Example}
\newtheorem{Remark}[Theorem]{Remark}
\newtheorem{Algorithm}[Theorem]{Algorithm}
\newtheorem{Procedure}[Theorem]{Procedure}
\newtheorem{History}[Theorem]{Historical Remark}
\newtheorem{Cicci}{Remark$\times$Cicci}
\newtheorem{teo}{Remark$\times$teo}
\renewcommand{\theCicci}{\Alph{Cicci}}
The aim of this note is to discuss  the following quite queer
\begin{Problem}\label{Prob} Given
\begin{itemize}
\item the   free non-commutative polynomial ring, ${\Cal P} := {\Bbb F}\langle X_1,\ldots,X_n\rangle$ {\em (public)},
\item a bilateral ideal ${\sf I}\subset {\Bbb F}\langle X_1,\ldots,X_n\rangle$ {\em (private)}, 
\item a finite set
$G := \{g_1,\ldots,g_l\}\subset{\sf I}$ of elements of the ideal  ${\sf I}$ {\em (public)},
\item a noetherian semigroup term-ordering $\prec,$ {\rm (private)}, on the word semigroup 
${\Cal T} := \langle X_1,\ldots,X_n\rangle$,
\end{itemize}
compute 
\begin{itemize}
\item[]  a finite subset $H\subset\Gamma({\sf I})$  of the Gr\"obner basis  $\Gamma({\sf I})$ 
  of ${\sf I}$ w.r.t. $\prec$ s.t., for each $g_i\in G$
 its {\em normal form} $NF(g_i,H)$ w.r.t. $H$ is zero,
\end{itemize}
 by means of a finite number of queries to
an oracle, which 
\begin{itemize}
\item[] given a term $\tau\in{\Cal T}$ returns its {\em canonical form} $\Can(\tau,{\sf I},\prec)$ w.r.t. the ideal ${\sf I}$ and the term-ordering $\prec$.
\qed\end{itemize}\end{Problem}

 
This queer problem has been suggested to us by  \cite{Bulygin} where a similar problem, but with stronger  assumptions, is faced in order to set up a  chosen-cyphertext attack  against the cryptographic system proposed in \cite{Rai}\footnote{Though we will breefly report on this application in Appendix we are not interested in dealing with it, preferring to  refer to the recent survey  \cite{Polly}.}.

The formulation of  Problem~\ref{Prob} is partially due to  the underlying application but is also due to the structure of the Gr\"obner bases in the non-commutative setting, which in general are infinite; however, even if we restrict to the noetherian setting of the (commutative) polynomial ring
${\Cal P} := {\Bbb F}[X_1,\ldots,X_n]$, we are unable (as we will show through easy counterexamples) to produce an algorithm which allows to return the (while finite) Gr\"obner basis of ${\sf I},$ unless we have some further informations   allowing to bound such basis; the best we can do is to solve the following reformulation:
 \begin{Problem}\label{Prob2} Given
\begin{itemize}
\item the  commutative polynomial ring, ${\Cal P} := {\Bbb F}[X_1,\ldots,X_n]$,
\item an ideal ${\sf I}\subset {\Bbb F}[X_1,\ldots,X_n]$, 
\item a noetherian semigroup term-ordering $\prec$ on the set of terms  
$${\Cal T} := \{ X_1^{a_1}\ldots X_n^{a_n}, (a_1,\ldots,a_n)\in{\Bbb N}^n\},$$  
\item a degree bound of the elements of the Gr\"obner basis $\Gamma(\sf I)$ of $\sf I$ w.r.t. $\prec$, {\em  i.e.} a  value $D\in{\Bbb N}$ satisfying $D\geq d(\sf I):= \max \{\deg (\gamma_i) : \gamma_i\in\Gamma(\sf I)\}$,
\end{itemize}
compute 
\begin{itemize}
\item  the Gr\"obner basis $\Gamma({\sf I})$ 
  of ${\sf I}$ w.r.t. $\prec,$ 
\end{itemize}
 by means of a finite number of queries to
an oracle, which 
 \begin{itemize}
\item given a term $\tau\in{\Cal T}$ returns its {\em canonical form} $\Can(\tau,{\sf I},\prec)$ w.r.t. the ideal ${\sf I}$ and the term-ordering $\prec$.
\qed\end{itemize}\end{Problem}

After recalling the basic notions and set up the notation (Section~\ref{Sec1}) we solve first Problem~\ref{Prob} (Section~\ref{Sec2}) and next Problem~\ref{Prob2} (Section~\ref{Sec3})  for which we propose a different, more combinatiorial, solution. 

We want to thank T. Moriarty and R.F. Ree for their precious apport.

\section{Notation and recalls on Gr\"obner Bases}\label{Sec1}

We consider a (non-necessarily commutative) monoid ${\Cal T}$ generated by the set of variables $\{X_1,\ldots,X_n\}$,
a field ${\Bbb F}$ and the monoid-ring 
${\Cal P} := \Span_{\Bbb F}({\Cal T})$.

\noindent For any set $F\subset{\Cal P}$ we denote ${\sf I}\subset{\Cal P}$ the (bilateral) ideal generated by $F$.

\noindent Each $f\in{\Cal P}$ can be uniquely
expressed as $$f = \sum_{\tau\in {\Cal T} } c(f,\tau) \tau\in\ {\Cal P};$$
and we call {\em support} of $f$ the set
$\supp(f)  := \{\tau\in {\Cal T}  : c(f,\tau) \neq 0\}.$

\noindent Moreover, fixing a  noetherian semigroup ordering  $\prec$ on ${\Cal T},$  the {\em leading term}, {\em leading coefficient} and  {\em leading monomial} of $f$ are ordinately:
$${\bf T}(f) := \max_\prec\{\tau\in \supp(f)\},\,
\, \lc(f) := c(f,{\bf T}(f)) \, \,  \,\hbox{and} \, \, \, {\bf M}(f) := \lc(f) {\bf T}(f).$$

For each ideal ${\sf I}\subset{\Cal P}$, we also consider
\begin{itemize}
\item  the  {\em semigroup ideal} ${\bf T}({\sf I}) := \{{\bf T}(f) : f\in{\sf I}\}$,  
\item the {\em Gr\"obner sous-escalier} ${\bf N}({\sf I}) :=  {\Cal T} \setminus {\bf T}({\sf I})$ , 
\item the vector-space
${\Bbb F}[{\bf N}({\sf I})] := \Span_{\Bbb F}({\bf N}({\sf I}))$,
\item ${\bf G}({\sf I})\subset {\bf T}({\sf I})$ the unique minimal basis of  ${\bf T}({\sf I}).$
\end{itemize}

We recall that  for  $f\in{\Cal P}$ and $G\subset{\Cal P},$
\begin{itemize}
\item  $f$ has  {\em Gr\"obner representation} 
in terms of $G$ 
if 
$$f = \sum_{i=1}^{\mu_f} c_i \lambda_i g_{j_i}\rho_i, \, \, c_i\in {\mathbb F}\setminus\{0\}, \lambda_i,\rho_i\in{\Cal T}, g_{j_i} \in G, \mu_f\in{\Bbb N}$$
with  
${\bf T}(f) = \lambda_1 {\bf T}(g_{j_1})\rho_1 \succ \cdots \succ \lambda_i {\bf T}(g_{j_i})\rho_i \succ \cdots.$
\item $h := NF(f,G,\prec)\in{\Cal P}$  is a {\em normal form}\index{Normal form} of $f$ w.r.t. $G$, if
\begin{itemize}
\item $f - h\in {\Bbb I}(G)$ has a  Gr\"obner representation in terms of $G$
and 
\item $h \neq 0 \then {\bf T}(h) \notin 
\left\{\lambda{\bf T}(g)\rho: 
 \lambda,\rho\in{\Cal T}, g \in G\right\} =: {\bf T}(G).$ 
\end{itemize}
\item For each $f\in{\Cal P},$ there is a unique {\em canonical form}
$$g := \Can(f,{\sf I},\prec)  = \sum_{t\in{\bf N}({\sf I})} \gamma(f,t) t 
\in  {\Bbb F}[{\bf N}({\sf I})]$$ 
s.t. $f - g \in {\sf I}$.
\item A Gr\"obner basis  of ${\sf I}$ is any set $\Gamma\subset{\sf I}$ s.t. $\{{\bf T}(\gamma) : \gamma\in \Gamma\}$ generates ${\bf T}({\sf I})$.
\item The {\em reduced Gr\"obner basis} of ${\sf I}$ is the set
$$\{\tau-\Can(\tau,{\sf I},\prec): \tau\in{\bf G}({\sf I})\}.$$

 \end{itemize}

\section{Oracle-supported Approximation of ${\Gamma}({\sf I})$}\label{Sec2}

Let us now specialize ${\Cal T}$ to be the word semigroup
${\Cal T} := \langle X_1,\ldots,X_n\rangle$
so that in particular the following holds:
\begin{itemize}
\item for each term $\upsilon\in{\Cal T}$ and  variables $X_l,X_r$ we have by definition
\begin{equation}\label{eq}
X_l\upsilon X_r\in{\bf G}({\sf I})  \iff X_l\upsilon\in{\bf N}({\sf I}), \upsilon X_r\in{\bf N}({\sf I}), X_l\upsilon X_r\in{\bf T}({\sf I});
\end{equation}
\item for each term $\upsilon\in{\Cal T}$ and each variable $X$ we have
\begin{equation}\label{eq2}
\omega = \upsilon X \in{\bf N}({\sf I})\implies  \upsilon\in{\bf N}({\sf I}), \, \,    
\omega = X\upsilon  \in{\bf N}({\sf I})\implies  \upsilon\in{\bf N}({\sf I}).
\end{equation}
 \end{itemize}

If we ask our oracle the value of $\Can(\tau,{\sf I},\prec)$\footnote{Or, in order to mask our question --- see the discussion on Bulygin assumption (B2) in the Appendix,
--- the values of $\Can(l_{\iota}\tau r_{\iota},{\sf I},\prec)$ where 
 $ l_{\iota},r_{\iota}\in{\Cal P}$  satisfy
$\tau = \sum_\iota l_{\iota}\tau r_{\iota}$, so that
$$\Can(\tau,{\sf I},\prec) = \sum_\iota \Can(l_{\iota}\tau r_{\iota},{\sf I},\prec).$$}
for any term $\tau\in{\Cal T}$, we can deduce whether 
\begin{enumerate}
\item $\tau\in{\bf T}({\sf I})$ in which case we obtain also $\Can(\tau,{\sf I},\prec)$, or
\item $\tau\in{\bf N}({\sf I})$ {\em i. e.} $\tau = \Can(\tau,{\sf I},\prec)$.
\end{enumerate}

\begin{Procedure}
We are assuming of having the sets $$\supp(g_j), g_j\in G,$$ so that, without needing to know the term-ordering $\prec$, we can deduce the sets $$T_j := \{\tau\in\supp(g_j) : \tau\nmid \omega, \forall \omega\in\supp(g_j)\}.$$

Since for each $j$, there are $\tau\in T_j, \lambda,\rho\in{\Cal T} : \tau = \lambda{\bf T}(f)\rho$ for some $f\in\Gamma({\sf I})$ {\em e.g.} $\tau := {\bf T}(g_j)\in{\bf T}({\sf I})$,  we can produce a scheme, based on Equation~(\ref{eq}), which in a finite number of steps produces an element of $\Gamma({\sf I})$; we choose the most suitable set $T_j$ then repeatedly we 
\begin{itemize}
\item pick an element $\tau\in T_j$, if $\tau\notin{\bf T}({\sf I}),$ simply remove it, otherwise:
\item for $\tau = X_l\omega\in{\bf T}({\sf I})$ we test whether $\omega\in{\bf T}({\sf I})$ in which case we set $\tau := \omega$ and repeat
 until we have an element $\tau =X_l\omega\in{\bf T}({\sf I})$ for which $\omega\in{\bf N}({\sf I})$;
\item now, for $\omega = \upsilon X_r\in{\bf N}({\sf I})$ we test whether $X_l\upsilon\in{\bf T}({\sf I}),$ in which case we set $\omega := \upsilon\in{\bf N}({\sf I})$ and repeat  until we have an element $X_l \upsilon X_r$ for which
$$X_l\upsilon\in{\bf N}({\sf I}), \upsilon X_r\in{\bf N}({\sf I}), X_l\upsilon X_r\in{\bf T}({\sf I})$$ {\em id est} 
$X_l\upsilon X_r\in{\bf G}({\sf I})$.
\end{itemize}

Remarking that we also have 
$${\bf G}({\sf I})\ni X_l\upsilon X_r \mid \tau\in\supp(g_j),$$
we can solve Problem~\ref{Prob} by a 
repeated application of the scheme above as follows:
set $H :=\emptyset$ and repeatedly
\begin{itemize}
\item  apply the scheme above thus obtaining an element $\tau\in{\bf G}({\sf I})$ and the polynomial
$\Can(\tau,{\sf I},\prec)$,
\item set $H := H\cup\{\tau-\Can(\tau,{\sf I},\prec)\}$, $G := \{NF(g,H) : g\in G\}$
\end{itemize}
until $G = \{0\}$.

At termination, which  is granted by noetherianity, the set $H$ satisfies the conditions required in Problem 1.
\end{Procedure}

Clearly, in the non-commutative case, where in general Gr\"obner bases are infinite, we can not hope to produce the whole basis of ${\sf I}.$ 

\section{Oracle-supported Deduction of ${\Gamma}({\sf I})$ (commutative case)}\label{Sec3}

We begin by observing that also in the commutative case ${\Cal P} = {\Bbb F}[X_1,\ldots,X_n],$ with $\deg(X_i)=1, \forall \, 1\leq i\leq n$, a strong solution returning the complete basis of an ideal ${\sf I}\subset{\Cal P}$ can not be produced, unless further knowledge is assumed: in fact, given  ${\sf I}\subset {\Bbb F}[X_1,\ldots,X_n]$ and a value $\delta\in{\Bbb N}, \delta< d(\sf I)$, in general there are smaller ideals (see Remark~\ref{teor})
${\sf J}\subsetneqq{\sf I}$ which satisfy
$$\{f\in{\sf I} : \deg(f) \leq\delta\} = \{f\in{\sf J} : \deg(f)\leq\delta\}.$$

  We recall the following definitions and facts: 
\begin{itemize}
\item For any $\tau\in{\Cal T}, 1\leq i\leq n$ the {\em $X_i$-th predecessor of} $\tau$ is $\frac{\tau}{X_i}$ if $X_i\mid\tau$, otherwise we say that $\tau$ does not have $X_i$-th predecessor.
\item
${\bf B}({\sf I})\subset {\bf T}({\sf I})$,  the  {\em border of the ideal}, is defined by 

${\bf B}({\sf I}):=\{\tau\in {\bf T} ({\sf I}) \, : \exists \, 1\leq i \leq n, \frac{\tau}{X_i}\in{\bf N}({\sf I}) \},$
\item ${\bf J}({\sf I})\subset {\bf T}({\sf I})$ the  {\em interior of the ideal},  is defined by 

${\bf J}({\sf I}):=\{\tau\in {\bf T} ({\sf I}) \, : \forall \, 1\leq i \leq n, \,  \frac{\tau}{X_i}\in{\bf T}({\sf I}) \},$ and 
\item  
the unique minimal basis of  ${\bf T}({\sf I}),$ ${\bf G}({\sf I})\subset {\bf B}({\sf I})$, is characterized as

${\bf G}({\sf I}):=\{\tau\in {\bf B} ({\sf I}) \, : \forall \, 1\leq i \leq n, \,  \frac{\tau}{X_i}\in{\bf N}({\sf I}) \}.$

 \item For each $f_1,f_2\in {\Cal P}$,
the {\em S-polynomial of $f_1$ and $f_2$}\index{S-polynomial}
is 
$$S(f_1,f_2) := 
\lc(f_2)^{-1} {\frac{\delta(f_1,f_2)}{{\bf T}(f_2)}} f_2 - \lc(f_1)^{-1}  {\frac{\delta(f_1,f_2)}{{\bf T}(f_1)}} f_1,$$
where $\delta := \delta(f_1,f_2) := \lcm({\bf T}(f_1),{\bf T}(f_2))$. 
\item A set $G = \{g_1,\ldots,g_s\}$ is a Gr\"obner basis of ${\Bbb I}(G)$ iff
for each $i<j$ the S-polynomial $S(g_i,g_j)$ has a Gr\"obner representation in terms of $G.$
\item (Buchberger's Second Criterion) 

For each $f,g,h\in{\Cal P} :{\bf T}(h) \mid  \lcm({\bf T}(f),{\bf T}(g))$,  if
both $S(f,h)$ and $S(g,h)$ have a Gr\"obner representation in terms of $G$, the same is true for
$S(f,g)$.
\item We also set $d(\sf I):=\max\{\deg(\zeta): \zeta\in G(\sf I).$
\end{itemize}

Let  then ${\sf J}\subset {\Bbb F}[X_1,\ldots,X_n] := {\Cal P}$ be an ideal,  $\prec$ a noetherian semigroup term-ordering, 
$\Gamma({\sf J}) = \{\gamma_1,\ldots,\gamma_s\}$ the Gr\"obner basis of ${\sf J}$ w.r.t. $\prec$ and
$\delta\in{\Bbb N}$ any degree value s.t.
$\delta \geq d({\sf J})+1$.

Enumerate  the variables and the Gr\"obner basis elements in such a way that
$X_1 \prec X_2 \prec \ldots \prec X_n$ and 
$$i < j \iff \mbox{ either }
\begin{cases} \deg(\gamma_i)>\deg(\gamma_j)\mbox{ or } \cr 
\deg(\gamma_i)=\deg(\gamma_j)\mbox{ and } {\bf T}(\gamma_i) \succ {\bf T}(\gamma_j).\end{cases}$$

Denoting
$$\Omega := \min_\prec\{\tau\in{\bf T}({\sf I}), \deg(\tau) = \delta+1\}$$ and $d_i :=\deg(\gamma_i) < \delta,$
we necessarily have
$$\Omega = X_1^{\delta+1-d_s}{\bf T}(\gamma_s).$$

 We also let $h_0 := \Omega-\Can(\Omega,{\sf J},\prec),$ so that $\lc(h_0)=1, {\bf T}(h_0)=\Omega = X_1^{\delta-d_s}{\bf T}(\gamma_s),$ and $h_i := X_2\gamma_i, 1\leq i \leq s.$ We obtain\footnote{Of course, our construction is indebted to the counterexample to Cardinal's Conjecture proposed in \cite{Mou}.}:
 
\begin{Proposition} With the above notation it holds
$H := \{h_0,h_1,\ldots,h_s\}$ is a  Gr\"obner basis w.r.t. $\prec$ of the ideal
${\Bbb I}(H) = X_2{\sf J}+(h_0)$.
\end{Proposition}

\pf Clearly if $S(\gamma_i,\gamma_j), 1\leq i<j\leq s,$ has the Gr\"obner representation in terms of $\Gamma({\sf J}), S(\gamma_i,\gamma_j)=\sum\limits_{\alpha=1}^{\mu_{ij}} c_\alpha \tau_\alpha \gamma_{\ell_{\alpha}},$ then  $S(h_i,h_j)=X_2\sum\limits_{\alpha=1}^{\mu_{ij}} c_\alpha \tau_\alpha \gamma_{\ell_{\alpha}}=\sum\limits_{\alpha=1}^{\mu_{ij}} c_\alpha \tau_\alpha h_{\ell_{\alpha}}$ is a Gr\"obner representation in terms of $H.$ 

\noindent Moreover, since $ \Omega={\bf T}(h_0) \, \hbox{and} \, {\bf T}(h_s)=X_2{\bf T}(\gamma_s) \mid \lcm({\bf T}(h_j),\Omega), 0\leq j\leq s$, as a direct consequence of Buchberger's Second Criterion, in order to prove the claim it is sufficient to show that the S-polynomial $S(h_s,h_0)$ between $h_0$ and $h_s$ has a Gr\"obner representation in terms of $H.$

By assumption there \, $\exists \, \mu=\mu_{h_0}, \alpha\in  {\Bbb N}, 1\leq\alpha \leq s,  c_\alpha\in{\Bbb F}\setminus\{0\}, \tau_\alpha\in{\Cal T},$ s.t. we have a  Gr\"obner representation
$$\sf J\ni h_0 =  \Omega-\Can(\Omega,\sf J,\prec) = \lc(\gamma_s)^{-1}  X_1^{D-d_s}\gamma_s + 
\sum_{\alpha=1}^\mu c_\alpha \tau_\alpha \gamma_{\ell_\alpha}$$
where $\gamma_{\ell_\alpha}\in \Gamma({\sf J})$ and
$$\Omega = X_1^{D-d_s}{\bf T}(\gamma_s) \succ \tau_1{\bf T}(\gamma_{\ell_1}) \succ \tau_2{\bf T}(\gamma_{\ell_2}) \succ \cdots;$$
 thus we trivially obtain the required Gr\"obner representation
\begin{eqnarray*}S(h_s,h_0) &=& \lc(h_0)^{-1} {\frac{\delta(h_s,h_0)}{{\bf T}(h_0)}} h_0 - \lc(h_s)^{-1}  {\frac{\delta(h_s,h_0)}{{\bf T}(h_s)}} h_s=\cr &=& X_2 h_0 -  \lc(\gamma_s)^{-1}  X_1^{D-d_s} (X_2 \gamma_s) \cr &=& X_2
\sum_{\alpha=1}^\mu c_\alpha \tau_\alpha \gamma_{\ell_\alpha} = \sum_{\alpha=1}^\mu c_\alpha \tau_\alpha h_{\ell_\alpha}.\end{eqnarray*}
\qed

\begin{Remark}\label{teor} For any ideal 
${\sf J}\subset{\Cal P}$,  noetherian semigroup term-ordering  $\prec$, and degree value
$\delta\in{\Bbb N}$  s.t.
$\delta \geq d(\sf J)+1$,
  the two ideals
${\sf I}_\delta := {\Bbb I}(H)$ and ${\sf I} := X_2{\sf J}$ 
 satisfy both:
$$\{f\in{\sf I}_\delta : \deg(f) \leq\delta\} = \{f\in{\sf I} : \deg(f) \leq\delta\} \, \, \textrm{and} \, \, \, {\sf I}\subset{\sf I}_\delta,$$
with 
$$d(\sf I_\delta)>\delta\geq d(\sf J)+1=d(\sf I).$$
Thus, the algorithm we are going to sketch below applied to the (unknown) ideal $\sf I_\delta$ returns the correct answer ${\sf I}_\delta$ if the input data satisfy $D\geq \delta+1$, but returns the wrong answer $\sf I$ if $\delta\geq D\geq d(\sf J)+1.$

\noindent That is, we actually need to assume to know an upper bound $D$ for $d(\sf I)$ and only deal with terms belonging to the {\em box} $$\Cal B(D):=\{X_1^{a_1}\cdots X_n^{a_n}\in\mathbb{T}: \, 0\leq a_i\leq D, \forall \, 1\leq i\leq n\}.$$\qed\end{Remark}

\medskip

We now give a combinatorial algorithm to solve Problem 2.

\bigskip

Let $\omega=X_1\cdot\ldots\cdot X_n,$ as $\omega^0 = 1\in{\bf  N}({\sf I}),$ we take iteratively $\omega^{i+1},i\in{\Bbb N},$ until either we find 
$j\in{\Bbb N},j\leq D$,  such that $\omega^{j-1}\in{\bf  N}({\sf I})$ and $\omega^j\in{\bf  T}({\sf I})$ or $\omega^{D}\in{\bf  N}({\sf I})$. In this last case we can deduce that ${\sf I}=(0)\footnote{In fact each term $\tau$ with $\deg(\tau)\leq D$ trivially satisfies $\tau \mid \omega^D$, i.e. $\omega^D\in{\bf  N}({\sf I})$ implies $ {\bf G}({\sf I})=\emptyset.$}$,  otherwise, for the found $j\in{\Bbb N}$ we begin deciding which of the following cases arises:
 \begin{enumerate}
 \item[Case 1]
  $\omega^j\in {\bf G}({\sf I})$ (i.e. all the predecessors of $\omega^j$ are in ${\bf  N}({\sf I})$),
  \item[Case 2] $\omega^j\in{\bf  B}({\sf I})\setminus{\bf G}({\sf I})$ (i.e. at most $n-1$ predecessors of $\omega^j$ are in ${\bf  N}({\sf I})$),
  \item[Case 3] $\omega^j\in{\bf  J}({\sf I})$ (i.e. all the predecessors of $\omega^j$ are in ${\bf  T}({\sf I})$).
  \end{enumerate}

To visualize the situation we identify ${\Cal T}$ with ${\Bbb N}^n$ thought as  $$\{\underline{x}=(x_1,\ldots,x_n)\in{\Bbb R}^n : x_i\in{\Bbb N}, 1\leq  i\leq n\};$$ by `line' (and one should better say `half-line') of ${\Cal T}$ we mean a set of aligned points of ${\Bbb N}^n\subset{\Bbb R}^n$ and similarly for  `plane', `hyperplane', `simplicial complex' etc..

We point out that :
\begin{itemize}
 \item[--] for $n=2$, ${\bf  B}({\sf I})$ is a `piecewise linear curve' ${\cal  C}({\sf I})$ consisting of contiguous horizontal and vertical `segments' from which all  the `convex' vertices are removed and possibly the leftmost vertical segment and the bottom horizontal one are `half-lines'\footnote{As ${\bf  B}({\sf I})\cup\{\rm all \, \, the \, \, convex \, \, vertices\}$ looks like the profile of a stair  A. Galligo introduced the term {\em escalier}.};
 \item[--] for $n\geq3$, $ {\bf  B}({\sf I})$ is a `simplicial complex'\footnote{Still called {\em escalier}.}, consisting of contiguous shares of `hyperplanes' each of them parallel to a `coordinate hyperplane' (the closest to a coordinate one  possibly being infinite) from which all  the `protruding' 
 $i$-th facets with $i\leq n-2$ are removed;
  \item[--] ${\bf  J}({\sf I})$ is the set of points lying above the {\em escalier}; 

 \item[--]  ${\bf G}({\sf I})$ consists of the `concave vertices' of the {\em escalier};
 \item[--]  ${\bf  N}({\sf I})$ is the set of points below the {\em escalier} (for this named {\em sous-escalier}).
 \end{itemize} 
 
 We will also call {\em $`0$-dimensional'}, $\ldots,$ {\em `$n-1$-dimensional'} {\em point of the {\em escalier}} a point lying on a vertex $,\ldots,$ on a $(n-1)$-facet (and not in a lower dimensional one) noticing that the elements of ${\bf G}({\sf I})$ are particular `$0$-dimensional' points.
 
From now on we will assume that 
$\exists j\in{\Bbb N}, j\leq D,$ such that $\omega^{j-1}\in{\bf  N}({\sf I})$ and $\omega^j\in{\bf  T}({\sf I}).$ 

\subsection{Two variables}

We distinguish between the three possible cases for $\omega^j:=X^jY^j$  
and, through several steps, we construct 
${\bf G}({\sf I}):$
\begin{enumerate}
\item[case 1] \, $\omega^j\in {\bf G}({\sf I})$ (the `line' $x=y$ meets ${\bf  T}({\sf I})$ in a `concave vertex' of the {\em escalier}),
\begin{enumerate}  
\item[] $I$ step: $t_1:=\omega^j=X^jY^j\in{\bf G}({\sf I})$ and we store it (it could be the only generator) 
\item[] $II$ step: 
starting from $t_1=\omega^j\in{\bf G}({\sf I})$ (found in step I), we need to consider $X^jY^{j+n}$ and  $X^{j+m}Y^j$ as $n,m\in{\Bbb N}^*$:
\item[a)] examine $X^jY^{j+n}$:
\begin{itemize}
\item [(i)] \, if $\forall \, n\leq D-j, \, X^{j-1}Y^{j+n}\in{\bf  N}({\sf I}),$ then 
there is no generator in ${\bf G}({\sf I})$
with $X$-exponent  $< j$; 
\item [(ii)] \, if $\exists \, \tilde{n}=\min \{n \, : \, 0<n\leq D-j,
X^{j-1}Y^{j+n}\in{\bf  T}({\sf I})\},$  we let $b_2:=j+\tilde{n}$ and 

- \, if $Y^{b_2}\in{\bf  T}({\sf I})$  then we set  $ \alpha_2 := 0$ 

- \, otherwise we set
$ \alpha_2 := \max\{\alpha\leq j-1 : \, X^{\alpha-1}Y^{b_2}\in{\bf  N}({\sf I})\},$ so that $t_{22}:=X^{\alpha_2}Y^{b_2}$, with $0\leq \alpha_2<j, b_2>j,$ is a new generator and we store it;
\end{itemize}

\item[b)] examine $X^{j+m}Y^j$:
\begin{itemize}
\item [(i)] \, if $\forall \, m\leq D-j, \, X^{j+m}Y^{j-1}\in{\bf  N}({\sf I}),$ then 
there is no generator in ${\bf G}({\sf I})$ with $Y$-exponent  $< j$; 
\item [(ii)] \, if $\exists \, \tilde{m}=\min \{0<m\leq D-j : \, X^{j+m}Y^{j-1}\in{\bf  T}({\sf I})\},$ 
we let $a_2:=j+\tilde{m}$ and 

- \, if $X^{a_2}\in{\bf  T}({\sf I})$  then we set  $ \beta_2 := 0$ 

- \, otherwise we set
$ \beta_2 := \max\{\beta\leq j-1 : \, X^{a_2}Y^{\beta-1}\in{\bf  N}({\sf I})\},$ so that $t_{21}:=X^{a_2}Y^{\beta_2}$, with $0\leq \beta_2<j, a_2>j$ is a new generator and we store it ;
\end{itemize}
$t_1$ is the only generator of ${\bf T}({\sf I})$ iff at  step $II$ hold both $a)(i)$ and $b)(i),$ otherwise  at least one further generator is found.
\end{enumerate}

\item[case 2] \,  $\omega^j\in{\bf  B}({\sf I})\setminus{\bf G}({\sf I}):$ have to distinguish whether the `line' $x=y$ meets ${\bf  T}({\sf I})$ in a `vertical' or `horizontal side' of the {\em escalier}: 
\begin{itemize}
\item [a)] $X^{j-1}Y^j\in{\bf  N}({\sf I}), X^jY^{j-1}\in{\bf T}({\sf I})$  
(`vertical side' case), 
\begin{itemize} 
\item[$I$ step]: - \, if $X^{j}\in{\bf  T}({\sf I})$  then we set  $\bar{\beta}_1 := 0$ 

 - \, otherwise we set
$$\bar{\beta}_1 := \max\{\beta< j : \, X^jY^{\beta-1}\in{\bf  N}({\sf I})\},$$ so that $\bar{t}_1:=X^jY^{\bar{\beta}_1}\in {\bf G}({\sf I})$ and we store it (possibly the only generator);

\item[$II$ step]: 
\begin{itemize}
\item[(j)] starting from   $\bar{t}_1:=X^jY^{\bar{\beta}_1}\in {\bf G}({\sf I})$, if $j<D$  we repeat the procedure described in case 1, step $II \, b) (i),(ii)$ possibly finding a new generator $\bar{t}_{21}:=X^{\bar{a}_2}Y^{\bar{\beta}_2}\in{\bf G}({\sf I})$ with $0\leq\bar{\beta}_2<\bar{\beta}_1<j, D\geq\bar{a}_2>j;$
\item[(jj)] starting from  $\omega^j$  we repeat the procedure described in case 1 step $II a) (i), (ii)$ possibly finding a new generator $\bar{t}_{22}:=X^{\bar{\alpha}_2}Y^{\bar{b}_2}\in{\bf G}({\sf I})$ with $0\leq\bar{\alpha}_2<j, D\geq\bar{b}_2>j;$
\end{itemize} 
\end{itemize}

\item [b)] $X^jY^{j-1}\in{\bf  N}({\sf I}), X^{j-1}Y^j\in{\bf T}({\sf I})$ 
(`horizontal side' case), 
\begin{itemize} 
\item[$I$ step]: -\, if $Y^{j}\in{\bf  T}({\sf I})$  then we set  $ \tilde{\alpha}_1 := 0$ 

\quad   - \, otherwise we set
$ \tilde{\alpha}_1 := \max\{\alpha<j : \, X^{\alpha-1}Y^j\in{\bf  N}({\sf I})\},$ so that $\tilde{t}_1:=X^{\tilde{\alpha}_1}Y^j\in{\bf G}({\sf I})$ and we store it (possibly the only generator);
\item[$II$ step]: 
\begin{itemize}
\item[(j)] starting from  $\tilde{t}_1:=X^{\tilde{\alpha}_1}Y^j
\in{\bf G}({\sf I})$, if $j<D$  we repeat the procedure described in case 1, step $II \, a)(i),(ii)$ possibly finding a new generator $\tilde{t}_{22}:=X^{\tilde{\alpha}_2}Y^{\tilde{b}_2}\in{\bf G}({\sf I})$ with $0\leq\tilde{\alpha}_2<\tilde{\alpha}_1<j, D\geq\tilde{b}_2>j;$
\item[(jj)] starting from  $\omega^j$  we repeat the procedure described in case 1 step $II b), (i), (ii)$ possibly finding a new generator $\tilde{t}_{21}:=X^{\tilde{\alpha}_2}Y^{\tilde{b}_2}\in{\bf G}({\sf I})$ with $0\leq\tilde{\alpha}_2<j, D\geq\tilde{b}_2>j;$
\end{itemize}
\end{itemize}
\end{itemize}
$\bar{t}_{1}$ (resp. $\tilde{t}_{1})$ is the only generator of ${\bf T}({\sf I})$ iff at  step $II \, a)$ (resp. $II \, b))$ hold both $a)(i)$ and $b)(i)$ of case 1 step $II$, otherwise  at least one further generator is added.

\item[case 3]  \, $\omega^j\in {\bf  J}({\sf I})$ (the `line' $x=y$ meets ${\bf  T}({\sf I})$ in a `convex vertex' of the {\em escalier}), \begin{itemize} 
\item[$I$ step]: 
by construction $\omega^{j-1}\in{\bf  N}({\sf I}),$ thus  
$X^{j-1}Y^j,X^jY^{j-1}\in{\bf  B}({\sf I})$ (the first one in a `horizontal' and the second one in a `vertical side' of the {\em escalier}), operating on them respectively like in case 2 $b)$ step $I$ and  case 2 $a)$ step $I,$ we get two generators:
\begin{itemize}
\item $\breve{t}_{12}:=X^{\breve{\alpha}_1}Y^j, \, 0\leq\breve{\alpha}_1<j,$
\item $\breve{t}_{11}:=X^jY^{\breve{\beta}_1}, \, 0\leq\breve{\beta}_1<j;$
\end{itemize}

\item[$II$ step]: 
\begin{itemize}
\item operating on  $\breve{t}_{12}$ like in case 1, step $II \, a) (i),(ii)$ we possibly find a new generator $\breve{t}_{22}:=X^{\breve{\alpha}_2}Y^{\breve{b}_2}$ with $0\leq\breve{\alpha}_2<\breve{\alpha}_1<j, D\geq\breve{b}_2>j$
\item operating on  $\breve{t}_{11}$ like in case 1, step $II \, b) (i),(ii)$ we possibly find a new generator $\breve{t}_{21}:=X^{\breve{a}_2}Y^{\breve{\beta}_2}$ with $0\leq\breve{\beta}_2<\breve{\beta}_1<j, 
D\geq\breve{a}_2>j;$
\end{itemize}
\end{itemize}

$\breve{t}_{11}$ and $\breve{t}_{12}$ are the only generators of ${\sf I}$ iff at step $II $ hold both $a) (i)$ and $b) (i)$ of case 1 step $II$, otherwise at least one further generator is added.

\item[all cases] $III$ and further steps

starting from the previous step generators  (all of type $t_{i2}:=X^{\alpha_i}Y^{b_i}$ with $0\leq\alpha_i<\ldots<j, D\geq b_i>\ldots>j$ or $t_{i1}:=X^{a_i}Y^{\beta_i}$ with $ 0\leq\beta_i<\ldots<j, 
D\geq a_i>\ldots>j)$ we operate like in 
case 2 step $II (j)$ while  $D> b_i$ and $D> a_i$

\noindent  The procedure stops because our possible degrees do not exceed the fixed bound $D$ and we don't miss any generator since we are following the {\em escalier} point by point.
 \end{enumerate}
 
 \begin{Example}\rm Let ${\Cal P} = {\Bbb F}[X,Y], \, \omega=XY.$ 
 \begin{enumerate}
 \item ${\sf I}=(X^2Y^2,XY^3,X^4Y,Y^8), D=8$. 
 
 \noindent We have $\omega^1\in{\bf  N}({\sf I}), \omega^2\in{\bf  T}({\sf I})$ and $XY^2, X^2Y\in{\bf  N}({\sf I}),$ thus $\omega^2\in {\bf G}({\sf I});$ considering $X^{2+m}Y, m\leq D-2$ and $XY^{2+n}, n\leq D-2$ we see that:  
 
 $\min\{n: \, XY^{2+n}\in{\bf  T}({\sf I})\}=1,$ with $Y^3, XY^2\in {\bf  N}({\sf I}),$ thus $XY^3\in {\bf G}({\sf I})$;  
 
 $\min\{m: \, X^{2+n}Y\in{\bf  T}({\sf I})\}=2,$ with $X^3Y, X^4\in {\bf  N}({\sf I})$ thus $X^4Y\in {\bf G}({\sf I})$.

\noindent Starting from $XY^3$ we see that  $\min\{n: Y^{3+n}\in {\bf  T}({\sf I})\}=5$
thus $Y^8\in{\bf G}({\sf I})$;
while,  starting from  $X^4Y$ we see that $X^{4+m}\in {\bf  N}({\sf I}), \forall \, m\leq D-4$, so that do not exist generators 
with null $Y$-exponent. 
 \item ${\sf I}=(X^3Y^2), D=5.$
  
 \noindent We have $\omega^1, \omega^2\in{\bf  N}({\sf I}), \, \omega^3\in{\bf  T}({\sf I})$ with $X^2Y^3 \in{\bf  N}({\sf I})$ and $X^3Y^2 \in{\bf  T}({\sf I})$  thus we have to consider $X^3Y^{3-q}, 0<q\leq 3$, as $X^3Y^{2}\in{\bf  B}({\sf I}), X^3Y\in{\bf  N}({\sf I})$  we have $X^3Y^2\in {\bf G}({\sf I});$ moreover as  $X^{3+m}Y\in{\bf  N}({\sf I}), \forall \,  m\leq D-3$ and $X^2Y^{2+n}\in{\bf  N}({\sf I}), \forall \,  n\leq D-2$ we have that $X^3Y^2$ is the unique generator.
 \item ${\sf I}=(X^2Y^4,X^4Y^3), D=7.$
 
 \noindent We have $\omega^1, \omega^2,\omega^3\in{\bf  N}({\sf I}), \, \omega^4\in{\bf  T}({\sf I})$ with $X^3Y^4,X^4Y^3 \in{\bf  B}({\sf I})$ 
 thus we have to consider $X^{4-p}Y^4,X^4Y^{4-q},  p,q\leq 4$, and we see that $X^4Y^{3}\in {\bf  G}({\sf I}), X^2Y^4\in {\bf  G}({\sf I})$ are the only generators of $\sf I$.
 \end{enumerate}
 \end{Example}

\subsection{$n\geq3$ variables}

Using the 2-variables case as a first inductive step, we consider $X_n$ as $n^{th}$ variable, added to $X_1,\ldots,X_{n-1}$. Assuming we are able to find all the minimal generators (up to the degree bound) of a monomial ideal in $n-1$ variables, we will  slice  ${\Cal T}$ in `hyperplanes' $x_n=j, j\leq D,$ and we will argue by considering the intersection $E_j$ of the {\em escalier} with each one of them. 
One of the following cases occurs: 
 \begin{itemize}
 \item $E_j$ has dimension  $i\leq n-2,$ so it does not contain any element of ${\bf G}({\sf I}),$
 \item $E_j$ is $n-1$-dimensional and so it contains some element of ${\bf G}({\sf I}),$ 
 \item $E_j=\emptyset.$
 \end{itemize}
 
\begin{Remark}\label{teo}  We point out explicitly that for any ${\sf I}\neq(0)$ there must exist at least one $j\in{\Bbb N}$ with $E_j$ hyperplanar.

Moreover, as we already remarked, $\omega^D\in{\bf N}({\sf I})\then {\sf I}=(0)$ and ${\bf N}({\sf I})=\emptyset.$ If, instead, for some $j\leq D, \omega^j\in{\bf T}({\sf I})$ then, necessarily, there is a $\tau\in{\bf G}({\sf I}), \tau\mid\omega^j$ and thus  $E_{j-h_1^-}\cap{\bf G}({\sf I})\neq\emptyset$ for some $h_1^-, 0 \leq h_1^- \leq j.$

It is however possible that for some $j\leq D$, $\omega^j\in{\bf T}({\sf I})$ and $E_{j+h}\cap{\bf G}({\sf I})=\emptyset$ for each $h, 0\leq h\leq D - j.$ This simply means that all generators of ${\bf T}({\sf I})$ have $X_n-$degree bounded by $j$ and that $E_j=E_{j+h}$ for each $h\in\mathbb{N}.$
\qed\end{Remark}
 \begin{enumerate}
 \item 
[Step I:] By applying the $n-1$-variables algorithm to $\omega^j$ (on the `hyperplane' $x_n=j)$ we find a set of terms ${\widetilde{\bf G}({\sf I})}_1$ from which, after cancelling all the terms $\sigma$ such that $\frac{\sigma}{X_n}\in{\bf  T}({\sf I}),$ we get  a set of terms ${\bf G}({\sf I})_{j\ldots j}$ for which
two possibilities arise:
\begin{itemize}
  \item[(i)] ${\bf G}({\sf I})_{j\ldots j}\neq\emptyset$ and we set ${\bf G}({\sf I})_1:={\bf G}({\sf I})_{j\ldots j},$ 
  \item[(ii)] otherwise, ${\bf G}({\sf I})_{j\ldots j}=\emptyset$ means that $E_j$ is $i\leq n-2$-dimensional
  and we have to iteratively consider $\omega_n^{+h}:=X_1^j\cdots X_{n-1}^jX_n^{j+h}, \forall \, \, h\leq D-j,$ 
and $\omega_n^{-h}:=X_1^j\cdots X_{n-1}^jX_n^{j-h}, \forall \, \, h\leq j,$ until we find  necessarily an  $E_{j-h}$ which is `hyperplanar' and possibly also an $E_{j+h},$ which is  `hyperplanar'; we then set\footnote{Notice that if ${\bf G}({\sf I})_{j\ldots j}\neq\emptyset$ we must think of $h_1^+=h_1^-=0$.}: 
\begin{itemize}
\item[--] $h_1^+:=\min \{h\leq D-j, E_{j+h}$ `hyperplanar'$\}$ (if it exists),
\item[--] $h_1^-:=\min \{h\leq j, E_{j-h}$ `hyperplanar'$\}$.
\end{itemize}
\end{itemize}

By applying the $n-1$-variables algorithm on both `hyperplanes'  $x_n=j+h_1^+$  and  $x_n=j-h_1^-$(noticing that by assumption $X_1^j\cdots X_{n-1}^jX_n^{j+h_1^+}$, $X_1^j\cdots X_{n-1}^jX_n^{j-h_1^-}\in{\bf  T}({\sf I}))$, after the above cancellation procedure, we get new sets of terms ${\bf G}({\sf I})_{j\ldots j}^{h_1^+} $ and ${\bf G}({\sf I})_{j\ldots j}^{h_1^-}$. As 
we observed in Remark \ref{teo} it can not happen $E_{j-h}\cap{\bf G}({\sf I})=\emptyset, \forall \, h\leq j,$ i.e. at least ${\bf G}({\sf I})_{j\ldots j}^{h_1^-}\neq\emptyset$
so that, setting : ${\bf G}({\sf I})_1^+:= {\bf G}({\sf I})_{j\ldots j}^{+h_1^+}$ and ${\bf G}({\sf I})_1^-:= 
{\bf G}({\sf I})_{j\ldots j}^{+h_1^-}\footnote{Of course if $\nexists \, h_1^+$ 
we set ${\bf G}({\sf I})_1^+:= \emptyset$ noticing that if ${\bf G}({\sf I})_1^+:= \emptyset$  do not exist generators with $X_n$-exponent $\geq j$. We also note that if ${\bf G}({\sf I})_{j\ldots j}\neq\emptyset$ we can think ${\bf G}({\sf I})_1={\bf G}({\sf I})_1^-.$}$, 
we get
$$\emptyset\neq{\bf G}({\sf I})_1:={\bf G}({\sf I})_1^+\cup{\bf G}({\sf I})_1^-,$$

\item[Step II]
  \begin{enumerate}
   \item[a)] 
  $\forall \, \, \sigma=X_1^{a_1}\cdots X_{n-1}^{a_{n-1}}X_n^{j-h_1^-}\in{\bf G}({\sf I})_1^-$ we move along the `line' $$\begin{cases} x_1-a_1=x_2-a_2\cr x_1-a_1=x_3-a_3\cr\vdots\cr x_n=j-h_1^--1,\end{cases},$$ with the following two possible issues: 
    \begin{itemize}
  \item[(i)] for all $ X_1^{a_1+l}\cdots X_{n-1}^{a_{n-1}+l}X_n^{j-h_1^-}\in{\bf G}({\sf I})_1^-$ and 
  $l\leq\max\{D-a_i\}$   
  it holds $$X_1^{a_1+l}\cdots X_{n-1}^{a_{n-1}+l}X_n^{j-h_1^--1}\in{\bf  N}({\sf I}),$$  that is the whole share of the `hyperplane' $x_n=j-h_1^-$ lying on  ${\bf  T}({\sf I})$ actually belongs to  ${\bf  B}({\sf I})$ (i.e. do not exist generators having $X_n$-exponent $< j-h_1^-).$
  \item[(ii)]  $\exists \, X_1^{a_1}\cdots X_{n-1}^{a_{n-1}}X_n^{j-h_1^-}\in{\bf G}({\sf I})_1^-$ and $$l_{a_1\ldots a_{n-1}}:=\min \left\{l\in{\Bbb N}^* \, : X_1^{a_1+l}\cdots X_{n-1}^{a_{n-1}+l}X_n^{j-h_1^--1}\in{\bf  T}({\sf I})\right\},$$ 
 that is the {\em escalier} does not exhaust ${\bf  T}({\sf I})\cap\{{\underline x}\in{\Bbb R^n} : x_n=j-h_1^-\}$ (i.e. some $X_1^{\alpha_1}\cdots X_{n-1}^{\alpha_{n-1}}X_n^{j-h_1^-}\in{\bf  J}({\sf I})$ and  do exist generators having $X_n$-exponent $< j-h_1^-).$ In this case we consider iteratively
$$X_1^{a_1+l_{a_1\ldots a_{n-1}}}\cdots X_{n-1}^{a_{n-1}+l_{a_1,\ldots,a_{n-1}}}X_n^{j-h_1^- -h}, h\leq j-h_1^-$$ 
until either we find
 $h_{a_1\ldots a_{n-1}}^-, 0< h_{a_1\ldots a_{n-1}}^- < j-h_1^-$
  with $$X_1^{a_1+l_{a_1\ldots a_{n-1}}}\cdots X_{n-1}^{a_{n-1}+l_{a_1\ldots a_{n-1}}}X_n^{j-h_1^--1-h_{a_1\ldots a_{n-1}}^-}\in{\bf  N}({\sf I})$$  
 (so that  $E_{j-h_1^--h_{a_1,\ldots,a_{n-1}}^-}$ is `hyperplanar'  thus containing  some generators of ${\sf I})$ or
 $X_1^{a_1+l_{a_1\ldots a_{n-1}}}\cdots X_{n-1}^{a_{n-1}+l_{a_1\ldots a_{n-1}}}\in{\bf  T}({\sf I})$ 
 in which case we set
 $h_{a_1\ldots a_{n-1}}^-=j-h_1^-$ (so that  $j-h_1^--h_{a_1\ldots a_{n-1}} = 0$ and still $E_0 = E_{j-h_1^--h_{a_1\ldots a_{n-1}}^-}$ is `hyperplanar'  thus containing  some generators of ${\sf I}).$
 
 We then set 
 $$h_2^-:=\min\limits_{X_1^{a_1}\cdots X_{n-1}^{a_{n-1}}X_n^{j-h_1^-}\in{\bf G}({\sf I})_1^-}\{h_{a_1\ldots a_{n-1}}^- \, \textrm{as above} \}.$$
By applying the $n-1$-variables algorithm on the `hyperplane' $x_n=j-h_1^--h_2^-$ (the nearest-below  which is $\parallel$ to $x_n=j-h_1^-$ and contains generators of ${\sf I}$) we find a set of terms ${\widetilde{\bf  G}}({\sf I})^{-h_2^-}$ from which we must erase all the terms whose $X_n$-predecessor lie in ${\bf  T}({\sf I}),$ getting, by construction, a non-empty:
 $${\bf G}({\sf I})^{-h_2^-}:={\widetilde{\bf  G}}({\sf I})^{-h_2^-}\setminus\{\sigma\in{\widetilde{\bf  G}}({\sf I})^{-h_2^-}\ \, : \frac{\sigma}{X_n}\in{\bf  T}({\sf I})\},$$
  which contains all the generators lying on the `hyperplane' $x_n=j- h_1^--h_2^- \, $ 
  \end{itemize}
and we let ${\bf G}({\sf I})_2^-:=\begin{cases} \emptyset \qquad \qquad \textrm{in case (i)}\cr{\bf G}({\sf I})^{-h_2^-} \quad \textrm{in case (ii)}\end{cases}.$ 
 
  \item[b)] If ${\bf G}({\sf I})_1^+\neq \emptyset$, we fix any $X_1^{a_1}\cdots X_{n-1}^{a_{n-1}}X_n^{j+h_1^+}\in{\bf G}({\sf I})_1^+:$ by iteratively applying (on each `hyperplane' $x_n=j+h_1^++h)$ the $n-1$-variables algorithm to $X_1^{a_1}\cdots X_{n-1}^{a_{n-1}}X_n^{j+h_1^++h}, j+h_1^++h\leq D$  we find a set of terms ${\widetilde{\bf G}({\sf I})}_2^{+h}$ from which, after cancelling all the terms $\sigma$ such that $\frac{\sigma}{X_n}\in{\bf  T}({\sf I}),$ we get a set ${\bf G}({\sf I})_2^{+h}$ and two  possibilities arise:
   
    \begin{itemize}
  \item[(i)] for all $h, j+h_1^++h\leq D, {\bf G}({\sf I})_2^{+h}=\emptyset$ which means that do not exist generators having $X_n$-exponent $>j+h_1^+;$
  \item[(ii)]  $\exists \, h^+_2=\min \{h, j+h_1^++h\leq D : {\bf G}({\sf I})_2^{+h}\neq\emptyset\}$ and ${\bf G}({\sf I})_2^{+h^+_2}$ gives all the generators contained in  the `hyperplane' $x_n=j+h_1^++h^+_2$ (the upper-nearest $\parallel$ to $x_n=j+h_1^+$ which contains generators).
 \end{itemize}

Then we let ${\bf G}({\sf I})_2^+:=\begin{cases} \emptyset \qquad \qquad \textrm{in case (i)}\cr{\bf G}({\sf I})^{+h_2^+} \quad \textrm{in case (ii)}\end{cases}$ 

We finally set ${\bf G}({\sf I})_2:={\bf G}({\sf I})^+_2\cup{\bf G}({\sf I})^-_2.$
 \end{enumerate}

  \item[Further Steps]: Starting from ${\bf G}({\sf I})_{i-1}= {\bf G}({\sf I})_{i-1}^+\cup{\bf G}({\sf I})_{i-1}^-$, $\forall \, \, i\geq3,$ we repeat:
  \begin{itemize}
    \item[--] if ${\bf G}({\sf I})_{i-1}^-\neq\emptyset$ for a fixed $\sigma \in {\bf G}({\sf I})_{i-1}^-$ all the procedures of  Step II a), possibly finding a non-empty ${\bf G}({\sf I})_{i}^-$ and the relative $X_n$-exponent $j-h_1^--\cdots-h_i^-$.
  \item[--]  if ${\bf G}({\sf I})_{i-1}^+\neq\emptyset$, for each $\sigma \in {\bf G}({\sf I})_{i-1}^+$ all the procedures of Step II b), possibly finding a non-empty ${\bf G}({\sf I})_{i}^+$. 
 \end{itemize}

The procedure stops because our possible degrees do not exceed the fixed bound $D\in {\Bbb N^*}$ that is we find an $n_D({\sf I})\in{\Bbb N}$ such that 
$${\bf G}({\sf I})_{\leq D}=\bigcup\limits_{i=1}^{n_D({\sf I})}{\bf G}({\sf I})_{i}$$ and we don't miss any generator since we have controlled the situation at each $x_n$-level.
 \end{enumerate}
 
   \begin{Example}\rm Let ${\Cal P} = {\Bbb F}[X,Y,Z], \, \omega=XYZ.$ 

\noindent ${\sf I}=(XY^3Z^4,Y^5Z^3,X^3Y^2Z^2,X^4Z), D=8$.
 
 \noindent We have $\omega^2\in{\bf  N}({\sf I}), \omega^3\in{\bf  T}({\sf I})$ with $X^3Y^3Z^2, X^3Y^2Z^3\in{\bf  T}({\sf I}),X^2Y^3Z^3\in{\bf  N}({\sf I}).$ 
 \begin{itemize}
 \item[Step I]
 We apply in the `plane' $z=3$ the 2-variables algorithm to $\omega^3=X^3Y^3(Z^3)$: as $X^2Y^3(Z^3)\in{\bf  N}({\sf I})$ and $X^3Y^2(Z^3)\in{\bf  T}({\sf I})$ we consider $X^3Y^{3-q}(Z^3)$, $q\leq 3$ until $X^3Y^{3-q}(Z^3)\in{\bf  B}({\sf I})$  and $X^3Y^{2-q}(Z^3)\in{\bf  N}({\sf I})$ or $q=3.$ Since $X^3Y^{2}(Z^3)\in{\bf  B}({\sf I})$  and $X^3Y(Z^3)\in{\bf  N}({\sf I})$ we take $X^3Y^{2}(Z^3)$ and we store it (recalling that $\omega^2\in{\bf  N}({\sf I})$). Starting from $X^3Y^{2}(Z^3)$ we consider $X^{3+m}Y(Z^3), m\leq 5,$ and, since $X^4Y(Z^3),X^4(Z^3)\in{\bf  B}({\sf I})$, we store $X^4(Z^3)$. Starting from $X^3Y^3(Z^3)$ we look whether $$\exists \, \nu:=\min\{n: \, X^2Y^{3+n}(Z^3)\in{\bf  T}({\sf I}), 3+n\leq 8\}$$ and we find $\nu=2$ as $X^2Y^5(Z^3)\in{\bf  B}({\sf I})$ from which, by considering  $X^{2-p}Y^5(Z^3), p\leq 2$ until $X^{2-p}Y^5(Z^3)\in{\bf  B}({\sf I})$ and $X^{1-p}Y^5(Z^3)\in{\bf  N}({\sf I})$ or $p=2,$ we obtain  
 $Y^5(Z^3)\in{\bf  T}({\sf I})$ and we store it. We stop here as the 2-variables algorithm on the `plane' $z=3$ does not produce other elements.
Dividing by $Z$ each $\sigma\in \{X^3Y^2Z^3,X^4Z^3,Y^5Z^3\}$ we get  ${\bf G}({\sf I})_1=\{Y^5Z^3\}$ (as $X^3Y^2Z^3,X^4Z^2\in{\bf  T}({\sf I})$).
 \item[Step II]
 \begin{itemize}
 \item[a)]
We look whether $\exists \, l_{0,5}:=\min\{l: \, X^{0+l}Y^{5+l}Z^2\in{\bf  T}({\sf I}), l\leq 8)\}$ and we get  $l_{0,5}=3$ (as $X^3Y^8Z^2\in {\bf  T}({\sf I})$ and $X^2Y^7Z^2\in {\bf  N}({\sf I}))$, we then consider $X^3Y^8(Z^2)$ on the `plane' $z=2$  and, by applying the 2-variables algorithm, we get $X^3Y^2Z^2\in {\bf  T}({\sf I})$ and $X^4Z^2\in {\bf  T}({\sf I})$ to be stored and, since dividing by $Z$, we get $X^3Y^2Z\in {\bf  N}({\sf I})$ while $X^4Z\in {\bf  T}({\sf I})$, we have
 ${\bf G}({\sf I})_2^-=\{X^3Y^2Z^2\}.$
\item[b)]  
 \noindent Let's now look to what happens on the `planes' $z=3+h, h\leq 5.$ Knowing that $X^3Y^3Z^4 \in{\bf  T}({\sf I})$ we must apply the 2-variables algorithm to $X^3Y^3(Z^4)$ on the `plane' $z=4$ obtaining as output the set 
 $$\{XY^3(Z^4),X^3Y^2(Z^4),X^4(Z^4)\}$$ and, as we have
 $X^3Y^2Z^3, X^4Z^3\in{\bf  T}({\sf I})$ but $XY^3Z^3\in{\bf  N}({\sf I})$ we set
 ${\bf G}({\sf I})_2^+=\{XY^3Z^4\}$ and finally ${\bf G}({\sf I})_2=\{XY^3Z^4$,$X^3Y^2Z^2\}$.
 \end{itemize}
  \item[Step III]
 \begin{itemize}
 \item[a)]
We look whether $\exists \, l_{3,2}:=\min\{l: \, X^{3+l}Y^{2+l}Z\in{\bf  T}({\sf I}),l\leq 6\}$ and we find  $l_{3,2}=1$ (as $X^4Y^3Z\in {\bf  T}({\sf I})$ and $X^3Y^2Z\in {\bf  N}({\sf I}))$, we then apply the 2-variables algorithm to $X^4Y^3(Z)$ on the `plane' $z=1$  finding only $X^4Z\in {\bf  B}({\sf I})$ to be stored and divided by $Z$ and, as $X^4\in {\bf  N}({\sf I}),$ we set
 ${\bf G}({\sf I})_3^-=\{X^4Z\}.$
\item[b)]  
Let's now look to what happens on the `planes' $z=4+h, h\leq 4,$ knowing that $XY^3Z^{4+h} \in{\bf  T}({\sf I})$ we apply the 2-variables algorithm to $XY^3(Z^{4+h}), h\leq 4;$ at each step we get
$$\{XY^3(Z^{4+h}),X^3Y^2(Z^{4+h}),X^4(Z^{4+h}), Y^5(Z^{4+h})\}$$
and since all elements are trivially to be discarded we get ${\bf G}({\sf I})_3^+=\emptyset$.
\end{itemize}
  \item[Further Steps]
Finally, since $X^{4+l}Y^{0+l}\in{\bf  N}({\sf I}), \forall \, l\leq 8,$ we deduce that there is no generator  with null $Z$-exponent, i.e. ${\bf G}({\sf I})_4^-=\emptyset$. Since we also have ${\bf G}({\sf I})_3^+=\emptyset,$ the algorithm terminates and we can conclude that 
${\bf G}({\sf I})=\{ XY^3Z^{4},X^3Y^2Z^{2},X^4 
Z,Y^5Z^3\}$.
 \end{itemize}
 \end{Example}


\section{A cryptographic application}

The survey \cite{Polly} reports on a class of cryptosystems whose scheme has been independently proposed by B.~Barkee {\em et  al.}  \cite{BCEMR} and by Fellows--Koblitz \cite{FK1,FK2,FK3,K}.
Such schemes are defined on the commutative polynomial ring ${\Cal P} = {\Bbb F}[X_1,\ldots,X_n]$ and
consist in:
\begin{enumerate}
\item writing down an easy-to-produce  Gr\"obner basis $\Gamma = \{\gamma_1,\ldots,\gamma_s\}$
generating an ideal ${\sf I} := {\Bbb I}(\Gamma)\subset {\Cal P}$ and
\item publishing a set $G := \{g_1,\ldots,g_l\}\subset {\sf I}$ of  polynomials  in ${\Cal P}$ and a set 
$$T := \{\tau_1,\ldots,\tau_m\}\subset{\bf  N}({\sf I}) = {\Cal T}\setminus{\bf  T}({\sf I})$$
of  {\em normal terms}
 belonging to the Gr\"obner 
{\em sous-escalier} of ${\sf I}$;
\item in order to send a message $M := \sum_{i=1}^m c_i\tau_i\in\Span_k(T)$, Bob (the sender) produces random polynomials $p_j\in {\Cal P}, 1\leq j\leq l, \deg(p_j) = \eth_j$, and encrypts $M$ as 
$C := M + \sum_{j=1}^l p_jg_j;$
\item  Alice (the receiver), possessing the Gr\"obner basis of ${\sf I},$ applies Buchberger's reduction to obtain
$\Can(C,{\sf I},\prec) = M =  \sum_{i=1}^m c_i\tau_i.$
\end{enumerate}

Rai \cite{Rai} proposed essentially the same system in the setting of the non-commuta\-tive polynomial ring  ${\Cal P} = {\Bbb F}\langle X_1,\ldots,X_n\rangle$: in his example the bilalteral ideal ${\sf  I}$  is  principal:
$${\sf I} := {\Bbb I}(\Gamma)\subset {\Cal P}, \, \, \Gamma = \{\gamma\}$$
and the published set $G := \{g_1,\ldots,g_l\}\subset {\sf I}$ is defined as
$g_i := h_i\gamma l_i$ for random elements $h_i,l_i\in{\Cal P}.$

We now describe a Bulygin-like (see \cite{Bulygin}) chosen-cyphertext attack on Barkee's  cryptosystems under  
the assumption of knowing
\begin{enumerate}
\renewcommand\theenumi{{\rm (B.\arabic{enumi})}}
\item the set ${\bf G}({\sf I}) := \{{\bf T}(\gamma_i) : \, 1 \leq i \leq s\}$ and
\item for each $\gamma_i\in \Gamma$, a set of pairs  $(s_i,t_i)$ of terms s.t. $s_iwt_i\notin {\bf T}({\sf I})$ for each $w\in\supp(\gamma_i).$
\end{enumerate}

Assuming the cryptoanalyst has temporary access to the decryption black box,
according Bulygin's attack, 
he then builds fake cyphertexts
$$C_i := s_{i}{\bf T}(\gamma_i)t_{i} + \sum_j p_j g_j q_j;$$
the decripted version of this message
being
$$\Can(C_i,{\sf I},\prec) = \Can(s_{i}{\bf T}(g_i)t_{i},{\sf I},\prec) = 
s_{i}\Can({\bf T}(g_i),{\sf I},\prec)t_{i}$$ thus  the attack allows him to read 
 $\gamma_i = {\bf T}(\gamma_i)-\Can({\bf T}(\gamma_i),{\sf I},\prec).$
 
 Before discussing the relation between Bulygin's assumption (B.1) and our oracle-based algorithm, let us consider the queer assumption  (B.2); it is justified by Bulygin as a tool for masking his attacks:
{\em Polynomial $t_{i},s_{i}$ are chosen for masking the ''fake'' cyphertext}  (\cite{Bulygin}, pg.2)

Assumption (B.2)  is  however completely useless: this ``masking'' in fact can be  performed simply by choosing  any set of polynomials $ l_{i\iota},r_{i\iota}\in{\Cal P}$ satisfying
${\bf T}(\gamma_i) = \sum_\iota l_{i\iota}{\bf T}(\gamma_i)r_{i\iota}$, thus we obtain
$$\Can({\bf T}(\gamma_i),{\sf I},\prec) = \sum_\iota \Can(l_{i\iota}{\bf T}(\gamma_i)r_{i\iota},{\sf I},\prec)$$
and we thus succeed in crashing the system via the  fake cyphertexts $l_{i\iota}{\bf T}(\gamma_i)r_{i\iota}$.


As regards assumption (B1), our investigation on the presented procedures was suggested by the aim of providing a tool to produce the set ${\bf G}({\sf I})$ and thus showing that assumption (B1) was unnecessary; however this is not true, except in the commutative case where we can cryptoanalyse a Barkee's scheme via our solution to Problem 2, provided we know a bound for the degrees.

{\em In fact we must stress that our solution of Problem~\ref{Prob} does not allow to reconstruct the set ${\bf G}({\sf I}),$ thus satysfying the necessary request (B1) by Bulygin, nor to cryptoanalyse a non-commutative Barkee's scheme:} all we can do is to produce a subset $H = \{h_1,\ldots,h_m\}\subset{\bf G}({\sf I})$ of the Gr\"obner basis $\Gamma({\sf I}) = \{\gamma_1,\ldots,\gamma_s\}$ --- used by Alice, via Buchberger's reduction, in order to read any message $M$ encrypted as $C = M+\sum_{j=1}^l p_j g_j q_j$ --- sufficient to produce a Gr\"obner representation 
$$g_j = \sum_i c_{ij} \lambda_{ij} h_{\iota_{ij}}\rho_{ij},
{\bf T}(g_j) =  \lambda_{1j} {\bf T}(h_{\iota_{1j}})\rho_{1j} \succ  \lambda_{2j} {\bf T}(h_{\iota_{2j}})\rho_{2j} \succ\ldots$$
of each  public element $g_j\in G$. Is this sufficient to obtain a  Gr\"obner representation of $C-M$? Of course no: in fact after we distribute the expression
$C-M=\sum_{j=1}^l p_j g_j q_j$ we obtain
$$C-M=\sum_{j=1}^L  \sum_i c_{j}  \lambda_{j} g_{\kappa_j}\rho_{j},
 \lambda_{j}, \rho_{j}\in{\Cal T}, c_j\in{\Bbb F}\setminus\{0\};$$ 
if we substitute each instance of $g_{\kappa_j}$ with its Gr\"obner representation deduced by our algorithm we simply have:
$$C-M=\sum_{j=1}^L  \sum_i c_{j} c_{i\kappa_j} \lambda_{j} \lambda_{i\kappa_j} h_{\iota_{i\kappa_j}}\rho_{i\kappa_j} \rho_{j};$$ 
thus if we properly reenumerate the summands we obtain a representation
$$C-M=\sum_{k=1}^K   d_k\lambda_k h_{\iota_k}\rho_{k}, \, \, \, 
\lambda_{1} {\bf T}(h_{\iota_1})\rho_{1} \succeq  \lambda_{2} {\bf T}(h_{\iota_2})\rho_{2} \succeq\ldots$$
but we
can not rule out equalities; thus we don't obtain
$${\bf T}(C-M) = \lambda_{1} {\bf T}(h_{\iota_1})\rho_{1} \succ  \lambda_{2} {\bf T}(h_{\iota_2})\rho_{2} 
\succ\ldots$$
and we cannot hope to successfully apply Buchberger reduction.

In fact, we can trivially build a theoretical counter-example by argueing as follows: assume that
$$\Omega:=\lambda_1{\bf T}(h_{\iota_1})\rho_{1} =\lambda_{2} {\bf T}(h_{\iota_2})\rho_{2}\succ\lambda_{3} {\bf T}(h_{\iota_3})\rho_{3} \,\, \,  \textrm{and} \, \, \, d_1\lc(h_{\iota_1})+d_2\lc(h_{\iota_2})= 0;$$
as a consequence, $l:=d_1\lambda_1h_{\iota_1}\rho_{1}+d_2\lambda_2h_{\iota_2}\rho_{2}\in{\sf I}$ necessarily satisfies ${\bf T}(l)\prec \Omega$ 
and has a Gr\"obner representation

$$l = \sum_{i=1}^I {\bar d_i}{\bar \lambda_i} \gamma_{\iota_i}{\bar \rho_i}, \qquad
{\bf T}(l) =  {\bar \lambda_1} {\bf T}(\gamma_{\iota_1}){\bar \rho_1}\succ\cdots$$
in terms of $\Gamma$ but not necessarily of $H$. Therefore, we can not discard the possibility that both 

\qquad ${\bar \lambda_1} {\bf T}(\gamma_{\iota_1}){\bar \rho_1}={\bf T}(l)\succ\lambda_{3} {\bf T}(h_{\iota_3})\rho_{3}$ and ${\bf T}(l)\notin\mathbb{I}({\bf T}(h):h\in H),$ 

\noindent so that $\gamma_{\iota_1}\notin H.$ In this unhappy, but realistic, situation we have the representation
$$C-M=\sum_{k=1}^K   d_k\lambda_k h_{\iota_k}\rho_{k}=l+ \sum_{k=3}^K   d_k\lambda_k h_{\iota_k}\rho_{k}=
\sum_{k=1}^I {\bar d_i}{\bar \lambda_i} \gamma_{\iota_i}{\bar \rho_i}+\sum_{k=3}^K   d_k\lambda_k h_{\iota_k}\rho_{k}$$
where 
$${\bar \lambda_1} {\bf T}(\gamma_{\iota_1}){\bar \rho_1}\succ{\bar \lambda_i} {\bf T}(\gamma_{\iota_i}){\bar \rho_i} \And {\bar \lambda_1} {\bf T}(\gamma_{\iota_1}){\bar \rho_1}\succ\lambda_{3} {\bf T}(h_{\iota_3})\rho_{3}\succeq\lambda_{k} {\bf T}(h_{\iota_k})\rho_{k}, \forall \, \, i,k,$$  so that necessarily ${\bf T}(C-M)={\bar \lambda_1} {\bf T}(\gamma_{\iota_1}){\bar \rho_1}\notin\mathbb{I}(\{{\bf T}(h):h\in H\})$ and we can not perform Bucheberger reduction.

On the other side, in the commutative case, each potential message $C$ necessarily satisfies
$$\deg(C)\leq \Delta:=\max\left\{\deg(\tau_i),\deg(g_j)+\eth, \tau_i\in T, g_j\in G\right\}$$
and thus $D := \Delta$ is a 'reasonable' guess  for degree bound $d({\sf I}).$
Of course the degree bound $ \Delta$  on the messages does not necessarily satisfy
$\Delta \geq d({\sf I}),$ so that our solution of Problem 2 would not cryptoanalyse Barkee's scheme using $D := \Delta$; however an implementation of Barkee's scheme in order to be protected  against it must assure $\Delta \ll d({\sf I}).$ 

While cryptoanalysing Barkee's schemes is an irrelevant task\footnote{Barkee's scheme was just a provocation aimed to address research towards sparse systems like the ones independently investigated, at the same time, by Fellows-Koblitz. As for their non-commutative generalizations, we simply wonder how it was possible that they have attracted attention, though an algorithm providing their cryptoanalysis [\cite{Pri}, Th. 13] was already available since 1996.} we would like to briefly point to a connected problem, which is equally irrelevant but at least is a combinatorial amusement. The technical tool used by the Barkee's scheme in order to {\em write down an easy-to-produce Gr\"obner basis} was later revealed in \cite{Mo} and simply consists into a combinatorial trick allowing, given any set of terms $\Upsilon:=\{v_1,\ldots,v_s\}\subset\Cal T,$ to produce a polynomial set $\Gamma:=\{\gamma_1,\ldots,\gamma_s\}$, satisfying ${\bf T}(\gamma_i)=v_i,$ and giving a Gr\"obner basis of the ideal it generates.

In principle, a Barkee's scheme could write down a term set $\Upsilon$ and the related easy-to-produce Gr\"obner basis $\Gamma,$ fix a value $D_0\ll d(\mathbb{I}(\Gamma))$, extract from $\Gamma$ the subset 

$\Gamma':=\{\gamma\in\Gamma: \, \deg(\gamma)\leq D_0\}$ with the corresponding term set

$\Upsilon':=\{{\bf T}(\gamma): \, \gamma\in\Gamma'\}=\{v\in\Upsilon: \, \deg(v)\leq D_0\}\subset\Upsilon$

\noindent and then produce the public set $G$ just using the elements belonging to $\Gamma'$ with
$$D_0<\Delta:=\max\{\deg(\tau_i),\deg(g_j)+\eth, \tau_i\in T, g_j\in G\}< d(\mathbb{I}(\Gamma)).$$
Recalling that our commutative procedure only deals with terms into the {\em box}
$$\Cal B(D):=\{X_1^{a_1}\cdots X_n^{a_n}\in\mathbb{T}: \, 0\leq a_i\leq D, \forall \, 1\leq i\leq n\},$$
and informally calling $D_0$-{\em badly-connected} a set of terms $\Upsilon$ such that, if we apply our procedure to it with the value $D := D_0 <\max\{\deg(v): \, v\in\Upsilon\}$ we are unable to produce the set $\Upsilon':=\{v\in\Upsilon: \deg(v)\leq D_0\}$, we remark that if $\Upsilon$ is $D_0$-badly connected, then in a Barkee's scheme, it would be nearly sufficient to make public a set $G\subset\mathbb{I}(\Gamma')$ in order to dwarf the use of our procedure in order to cryptoanalyse it.

The question, then, becomes the existence of badly connected sets of terms; we have the strong impression that the answer is negative\footnote{Consider the 2-variable case;  in a minimal Gr\"obner basis $\Gamma$, for any two elements $X^aY^b, X^cY^d\in{\bf G}(\Gamma)$ 
$a<c$ imples $b>d$.

Thus, if $X^aY^b, X^cY^d\in{\bf G}(\Gamma)$ are $D_0$-badly connected, there must be an element $X^eY^f\notin\Cal B(D)$
and which satisfyes $a<e<c, b>f>d$.
For such elements necessarily either $D_0 < e < a$ or $D_0 < f < b$, contradicting the assumption that 
$\deg(X^aY^b), \deg(X^cY^d) \leq D_0$.}. Nevertheless, as we said above, we consider irrelevant to devote some time to this task.


\end{document}